\documentclass[a4paper]{article}
\usepackage{Preamble}

\title{\LARGE \bf An Adjoint-based Numerical Method for a class of nonlinear Fokker-Planck Equations}

\author{Adriano Festa$^{1}$, Diogo A. Gomes$^{2}$ and Roberto M. Velho$^{2}$
\thanks{*D. Gomes was partially supported by KAUST baseline and start-up funds and 
KAUST SRI, Uncertainty Quantification Center in Computational Science and Engineering.  
        A. Festa was partially supported by the Haute-Normandie Regional Council via the M2NUM project and by ANR HJNet (ANR-12-BS01-0008-01)}
        \thanks{$^{1}$ Institut National de Sciences Appliqu\'ees,  LMI lab
Avenue de l'Universit\'e, 76800 Saint-\'Etienne-du-Rouvray,    Rouen, France
        {\tt\small adriano.festa@insa-rouen.fr}}
\thanks{$^{2}$King Abdullah University of Science and Technology (KAUST), CEMSE Division, Thuwal 23955-6900. Saudi Arabia, and  
        KAUST SRI, Center for Uncertainty Quantification in Computational Science and Engineering.
        {\tt\small diogo.gomes@kaust.edu.sa}, {\tt\small roberto.velho@gmail.com}}
}

\begin{document}
\maketitle
\thispagestyle{empty}
\pagestyle{empty}

\begin{abstract}
Here, we introduce a numerical approach for a class of Fokker-Planck (FP) equations. These equations are the adjoint of the linearization of Hamilton-Jacobi (HJ) equations.  
Using this structure, 
we show how to transfer the properties of schemes for HJ equations to the FP equations. Hence, we get numerical schemes with desirable features such as positivity and mass-preservation.
We illustrate this approach in examples that include mean-field games and a crowd motion model.
\end{abstract}


\section{Introduction}
Fokker-Plank (FP)  equations model the time evolution of a probability density. 
The general set up is as follows. 
Given an open subset of $\RR^d$, $\Omega$, a terminal time, $T>0$,
and a (\emph{drift})
 vector field, $b(x,t):\Omega \times[0,T]\rightarrow \Omega$, 
 we seek to find a time-dependent probability distribution,  
 $\rho:\Omega \times [0,T] \rightarrow \RR$, solving 
\begin{equation}\label{nonlinearFP}
\left\{\begin{array}{ll}
\partial_t \rho -\varepsilon\Delta \rho+ \mbox{div}(b(x,t) \rho) =0  &\mbox{in } \Omega \times [0,T], \\[6pt]
\rho(\cdot, 0)= \rho_0(\cdot) &\mbox{in }  \Omega.
\end{array}\right.
\end{equation}
In addition, we supplement the above problem with 
 boundary conditions on $\partial \Omega\times [0,T]$, where $\partial \Omega$ is 
 the boundary of $\Omega$. 

The Fokker-Planck equation was introduced in  statistical mechanics. Yet, 
this equation has multiple applications in economics \cite{aime10,Gueant09}, crowd motion models \cite{hughes2000flow,Lachapelle10}, and biological models \cite{chavanis2008nonlinear, goudon1998fokker}. Due to the complex structure of those equations, the computation of explicit solutions is not possible. Hence, effective numerical methods to approximate solutions of FP equations have a broad interest.

Here, we propose a technique to obtain approximation schemes for FP equations using their representation as the adjoint of the linearization of Hamilton-Jacobi (HJ) equations. In this way, all monotone numerical schemes proposed in the context of HJ equations give rise to consistent schemes for FP equations. In particular, these schemes preserve positivity and total mass, as required by the nature of the problem.

 Previously, the adjoint structure of the FP equation was used by several authors, for example, in \cite{AchdouCapuzzo10} and in \cite{AchdouCapuzzoCamilli12}. In those references, the authors propose a finite-difference scheme which is the adjoint of the linearization of the upwind scheme used to approximate a convex Hamiltonian. In \cite{carlinisilvafp2016} and \cite{carlini2016DGA}, the authors propose a semi-Lagrangian numerical method using a slightly different procedure, but based on a similar principle.
 
 The main contribution of the present paper is to show how to use the adjoint structure with a wide class of numerical solvers, and without limitations on the problem dimension. 
 Further, in contrast to the above references, we do not discretize the time variable. Thus, the evolution
 in time corresponds to a system of ordinary differential equations (ODE). These can be solved with different methods, depending on the smoothness of the solution and desired accuracy. Finally, 
 the implementation of our method uses a symbolic-numeric approach. Here, the numerical schemes are created by exact formula manipulation, thus reducing the implementation time and complexity.  
 
\smallskip
\noindent {\bf Outline of the paper.} We end this introduction with an outline of this paper. The adjoint structure is examined in Section \ref{eaa}. Next, in Section \ref{eee}, we proof key features of the method:  positivity and mass-conservation. In Section \ref{num}, we describe the numerical method and its properties. Some sample schemes are studied in detail.
Finally, in Section \ref{mf}, consider some problems where our schemes apply. 
These included mean-field games and a crowd motion model.

\section{Adjoint structure}\label{eaa}
The relation between a FP equation and its adjoint equation is well known. In recent works, \cite{evans2010adjoint, gomes2016regularity,MR3303946,MR3146863,MR3092361,MR2873233,MR2796233}, this relation was used to study regularity properties, vanishing viscosity limits, and rates of convergence of numerical methods.  
Those results are based on the observation that a FP equation is the adjoint of the linearization of a certain HJ equation.


\subsection{Adjoint structure}
Here, we discuss
the relation between FP and HJ equations. First, we consider the HJ operator
\begin{equation}\label{hj1}
HJ(u) := - u_t(x,t) + H(x,Du(x,t)) - \varepsilon \Delta u(x,t),
\end{equation}
with the Hamiltonian $H=H(x,p):\RR^d\times \RR^d\rightarrow~\RR$. Further, we define the nonlinear generator $$A^{HJ} u :=  H(x,Du(x,t)) - \varepsilon \Delta u(x,t).$$ Here, we write $Du = D_x u$ for the gradient in the variable $x = (x_1,\cdots,x_d)$. The parameter $\varepsilon$ is called the viscosity.

To linearize \eqref{hj1} around $u_0$, we expand $u = u_0 + \lambda w$,  then take the derivative in $\lambda$, and, finally, consider the limit $\lambda \to 0$.
For now, we proceed formally to compute this linearization. Later, we discuss
functional spaces and boundary conditions. 

The expansion $HJ(u_0 + \lambda w)$ gives
\begin{multline*}
- \partial_t (u_0 + \lambda w) + H(x,D(u_0 + \lambda w)) - \varepsilon \Delta(u_0 + \lambda w)\\[6pt]
= - {(u_0)}_t  - \lambda w_t + H(x,Du_0 + \lambda Dw) - \varepsilon \Delta u_0 - \lambda \varepsilon \Delta w.
\end{multline*}
By taking the derivative of the preceding expression with respect to $\lambda$, and letting $\lambda \to 0$, we obtain the operator
\begin{equation}\label{linearized_eq}
L(w):= - w_t + D_p H(x,Du) \cdot Dw - \varepsilon \Delta w,
\end{equation}
the linearization of the HJ operator. The (linear) generator of $L$ is
$$A^L w := D_pH(x,Du) \cdot Dw - \varepsilon \Delta w.$$

Finally, we compute
the adjoint of $L$ by integration by parts. We fix smooth functions, $w$ and $\rho$, and
derive the identity
\begin{align}\label{aeq}
& \iint\limits_{[0,T] \times \Omega} (-w_t + D_pH(x, Du) \cdot Dw - \varepsilon \Delta w) \ \rho \\\notag
= & \iint\limits_{[0,T] \times\Omega} \left(\rho_t - \div_x(D_p H(x, Du) \ \rho)
 - \varepsilon \Delta \rho \right) w\\\notag
+ & \iint\limits_{[0,T] \times \partial \Omega}\left(D_p H(x, Du) \right)\cdot n \ \rho \ w + \varepsilon \frac{\partial \rho}{\partial n} w - \varepsilon \rho \frac{\partial w}{\partial n}\\\notag
- & \int\limits_{\Omega} \rho(x,T) \ w(x,T) - \rho(x,0) \ w(x,0),
\end{align}
where  $n$ is the normal  vector to the boundary, $\partial \Omega$. The last calculation shows  that the adjoint of $L$ is the following FP operator
\begin{equation}\label{fp1}
L^* \rho := \rho_t - \div_x(D_p H(x, Du) \ \rho) - \varepsilon  \Delta \rho,
\end{equation}
whose generator is $A^{FP} \rho := - \div_x(D_p H(x, Du) \ \rho) - \varepsilon  \Delta \rho$.

\subsection{Boundary conditions}
Now, we address the boundary conditions for \eqref{nonlinearFP} on $\partial \Omega\times [0,T]$. 
The discussion of  initial conditions is straightforward. 
Two common boundary conditions for FP equations are Dirichlet data
and a prescribed flow via Neumann conditions. 
Typically, the Dirichlet data vanishes on the boundary. These boundary conditions correspond to the case where particles exit once they reach the boundary.
The prescribed flow case represents a current of particles or agents crossing the boundary. 
With a zero flow, the mass is conserved. 

Each of these choices
of boundary conditions
determines cancellations in the boundary integrals 
in \eqref{aeq}. This suggests different functional spaces for the HJ operator, its linearized version,  and its adjoint, the FP operator.

The first case corresponds to a FP equation with Dirichlet boundary conditions:
\begin{equation*}
\left\{\begin{array}{ll}
\rho_t(x,t) - \div(D_p H(x, Du) \ \rho) = \varepsilon  \Delta \rho, \hspace{0.15cm} \mbox{in } \Omega \times [0,T],& \\[6pt]
\rho(\cdot,t) = 0, \hspace{3.94cm} \mbox{on } \partial \Omega \times [0,T].&  
\end{array}\right.
\end{equation*}
We consider the HJ operator on a functional space with the boundary conditions
\begin{equation*}
\left\{\begin{array}{ll}
- u_t(x,t) + H(x,Du(x,t)) - \varepsilon \Delta u(x,t), \hspace{0.15cm}  \mbox{in } \Omega \times [0,T],& \\[6pt]
u(\cdot,t) = g_1(\cdot,t), \ \ \text{for any } g_1, \ \hspace{1.39cm} \mbox{on } \partial \Omega \times [0,T],& 
\end{array}\right.
\end{equation*}
and the linearized operator as
\begin{equation*}
\left\{\begin{array}{ll}
- w_t + D_p H(x,Du) \cdot Dw - \varepsilon \Delta w, \hspace{0.15cm} \mbox{in } \Omega \times [0,T],& \\[6pt]
w(\cdot,t) = 0, \hspace{3.25cm} \mbox{on } \partial \Omega \times [0,T].& 
\end{array}\right.
\end{equation*}

The second case corresponds to a FP equation with a flux through the boundary
\begin{equation*}
\left\{\begin{array}{ll}
\rho_t(x,t) - \div(D_p H(x, Du) \ \rho) - \varepsilon  \Delta \rho(x,t), \hspace{0.15cm} \mbox{in } \Omega \times [0,T],& \\[6pt]
D_p H(x, Du) \ \rho+\varepsilon
\frac{\partial \rho}{\partial n}(x,t) = g_2(x,t), \hspace{0.65cm} \mbox{on } \partial \Omega \times [0,T],&  
\end{array}\right.
\end{equation*}
where $g_2$ is the desired in/out-flow through $\partial \Omega$. We can consider diverse boundary conditions for the HJ operator: Dirichlet type, state-constraint, reflection at the boundary, and Neumann type.
In the following example, we use Neumann conditions. The Hamilton-Jacobi operator is
\begin{equation*}
\left\{\begin{array}{ll}
- u_t(x,t) + H(x,Du(x,t)) - \varepsilon \Delta u(x,t), \hspace{0.15cm} \mbox{in } \Omega \times [0,T],& \\[6pt]
\frac{\partial u}{\partial n}(x,t) = 0, \hspace{3.88cm} \mbox{on } \partial \Omega \times [0,T],& 
\end{array}\right.
\end{equation*}
with the corresponding linearization
\begin{equation*}
\left\{\begin{array}{ll}
- w_t + D_p H(x,Du) \cdot Dw - \varepsilon \Delta w, \hspace{0.15cm} \mbox{in } \Omega \times [0,T],& \\[6pt]
\frac{\partial w}{\partial n}(\cdot,t) = 0, \hspace{3.075cm} \mbox{on } \partial \Omega \times [0,T].& 
\end{array}\right.
\end{equation*}

We do not address the initial conditions for the above operators because we use them only to discretize in space the HJ generator.

A nonlinear FP equation is related to the solution of a stochastic differential equation of McKean-Vlasov type (or mean-field type), see \cite{MR0221595,MR0233437,MR1431299,MR1108185}. More precisely, we consider the stochastic differential equation (SDE)
\begin{equation}\label{Mackeanvlasov}
\left\{\begin{array}{ll}
d X(t)= b(X(t), \rho(X(t),t),t)\,d t+ \sqrt{2\varepsilon}\, d W(t),& \\[6pt]
X(0)= X^0,&
\end{array}\right.
\end{equation}
where $b: \RR^d \times \RR_+ \times \RR_+ \to \RR^d$ is a regular vector-valued function,  $X^0$ is a random vector in $\RR^d$, independent of the Brownian motion $W(\cdot)$, with density $\rho_0$, and  $\rho(\cdot,t)$ is the density of $X(t)$. It can be shown  (see \cite{MR1653393}) that under certain growth conditions for $b$ \eqref{Mackeanvlasov} admits a unique solution and $\rho$ is the unique classical solution of the nonlinear FP equation
\begin{equation*}\label{nonlinearFP2}
\partial_t \rho -\varepsilon\Delta \rho+ \mbox{div}(b(x,\rho,t) \rho) =0.
\end{equation*}
Therefore, if we set $b(x,m,t):=-D_p H(x,Du)$ and impose appropriate boundary conditions,  \eqref{Mackeanvlasov} provides a probabilistic interpretation of the optimal trajectories for \eqref{hj1}. 
With Dirichlet conditions, those trajectories end at the boundary; for zero-flux conditions, they are reflected, see \cite{Bossy2004}, and \cite{Gobet2001}.

%
%

\begin{remark}
	Our methods can be extended to study stationary FP equations. In this case, the associated  Hamilton-Jacobi operator is stationary. Small modifications can be added to the HJ operator to guarantee the existence of solutions.
\end{remark}

\section{Properties}\label{eee}
In this section, we show that the evolution of an initial density by the FP equation preserves positivity and mass. We use arguments from nonlinear semigroup theory to illustrate how these properties are related to corresponding
properties of the Hamilton-Jacobi equation. The arguments detailed here
are valid without any substantial changes for the discretized problems. 

We denote by 
 $\langle f, g \rangle = \int_{\Omega} f \ g$ the duality product, and by $S_t$ the semigroup associated to the linearized operator \eqref{linearized_eq}. This semigroup preserves order and commutes with constants. 
We define the adjoint $S^*_t$  by
\[
\langle S^*_tu, v\rangle= \langle u, S_tv\rangle.
\]

We have then the following results:
\begin{pro}[Positivity]
The evolution of the initial density $\rho_0$ through the adjoint semigroup, $S_t^*$, preserves positivity. Denote by $w_T$ the terminal condition for the linearized operator. Then, for $w_T \geq 0$, and $\rho_0 \geq 0$; we have $L^*_t \rho \geq 0$, for all $t \in [0,T]$.
\end{pro}

\noindent \emph{Proof: } First, note that $w_T \geq 0$ implies $S_t  w_T \geq 0$. This follows from the maximum principle for HJ equations. Thus, for $w_T \geq 0$, we have $$ \langle S_t^* \rho, w_T \rangle = \langle \rho, S_t w_T \rangle \geq 0,$$
since $\rho \geq 0$, and $S_t w_T \geq 0$. Accordingly, $S_t^* \rho \geq 0$.
\cqd

\begin{pro}[Conservation of Mass]
Let $\rho_0$ be the initial probability distribution, i.e. $\int_{\Omega} \rho_0 = 1$. Then, for all $t \in [0,T]$, the evolution of this probability measure through the adjoint semigroup, $S_t^* \rho_0$, is also a probability measure.
\end{pro}

\noindent \emph{Proof: }  First, observe that $S_t 1 = 1$. Then,
\begin{align*}
\int\limits_{\Omega} S_t^* \rho_0 = \langle S_t^* \rho_0,1 \rangle = \langle \rho_0, S_t 1 \rangle = \langle \rho_0,1 \rangle = \int\limits_{\Omega} \rho_0 = 1.
\end{align*}
\cqd

We conclude this section with some remarks.
\begin{remark}
	 In the computations of the previous sections, we assume that $Du(x,t)$ does not depend on $\rho$. Further, the relation between a general FP equation whose drift depends on the density, and its associated HJ equation is still a research topic. We do not address this case in the present work. Still, particular cases of drift depending on the density and numeric approaches to solve them are discussed in the literature, see for instance \cite{MR2724518}, and \cite{MR2966923}.
\end{remark}

\begin{remark}
	 If the viscosity vanishes (~$\varepsilon~=~0$~), the same approach is valid. A first-order HJ operator gives rise to a  continuity equation (CE), i.e. a FP equation without viscosity. This case is considered in section \ref{mf}, where we extend our numerical scheme to address systems of partial differential equations (PDEs). Those systems arise in multiple applications such as mean-field games (MFG), population models, traffic flow problems, and modeling in chemotaxis.
\end{remark}

\section{Numerical Approach}\label{num}

Our numerical approach relies on the relation between the HJ framework and the corresponding adjoint FP equation. Given a semi-discrete (discrete in space) numerical scheme for \eqref{hj1}, the same scheme can be used to construct a consistent approximation for \eqref{fp1}.

 Before proceeding, we define additional notation. To simplify, we consider a scheme for the case where the domain $\Omega$ is $\mathbb T^2$ (2-D torus). Let $\mathbb T^2_{\Delta x}$ be an uniform grid on $\mathbb T^2$ with constant discretization parameter $\Delta x>0$. Let $x_{i,j}$ denote a generic point in $\mathbb T^2_{\Delta x}$.  The space of  grid functions defined on  $\mathbb T^2_{\Delta x}$ is denoted by $\mathcal{G}(\mathbb T^2_{\Delta x})$, and the functions $U$ and $M$ in $\mathcal{G}(\mathbb T^2_{\Delta x})$ (approximations of respectively $u$ and $\rho$)  
are called $U_{i,j}$ and $M_{i,j}$, when evaluated at $x_{i,j}$.

 We utilize a semi-discrete numerical scheme $N(x,p):\mathbb T^2_{\Delta x}\times \RR^d\rightarrow \mathbb R$ monotone and consistent to approximate the operator $H(x,p)$, such that $U$ is the solution of the ODE  
\begin{equation}
U_t=N(x,\mathcal D U),
\end{equation}
where $\mathcal D U$ is a discretization of the gradient operator on $U$.
Thanks to the adjoint structure, we modify this scheme to approximate the solution of \eqref{fp1}. The discrete approximation, $M$, is the solution of the following system of ODE $$ M_t = K,$$ where
\begin{equation}\label{fp3}
K(x,\mathcal D U,M):=(D_p N(x,\mathcal D U))^T M+\varepsilon \Delta_d M.
\end{equation}
Here, the nonlinear part of the operator corresponds to
the discrete operator $D_p N(x_{i,j},\mathcal D U)$;  $\Delta_d M$ is a discretization of the Laplacian. We note that this operator depends on the monotone approximation scheme used to discretize the HJ equation, and can be computed numerically or using a symbolic differentiation operator. This is the case in our examples in section \ref{mf}.

We stress that the features of positivity and mass conservation are valid at the discrete level. This is a consequence of the semigroup arguments in section~\ref{eee}, independently of the manner the space or time are discretized.

%

\subsection{Finite Differences}\label{ssct:FinDif}
Now, we consider an explicit scheme using our method. We describe  an upwind discretization for the Hamiltonian, which we assume to be
\begin{equation}\label{ham}
H(x,p)=g(x)+|p|^\alpha, \quad \alpha>1.
\end{equation}
We define the standard finite-difference operators as
$$(\mathcal D_1^\pm u)_{i,j}=\frac{u_{i\pm 1,j}-u_{i,j}}{\Delta x},  \hbox{ and } (\mathcal D_2^\pm u)_{i,j}=\frac{u_{i,j\pm 1}-u_{i,j}}{\Delta x},$$
and 
$$\Delta_d u=\frac{1}{\Delta x^2}\left( 4 u_{i,j}-u_{i+1,j}-u_{i,j+1}-u_{i-1,j}-u_{i,j-1} \right).$$

 The  approximation of the operator $H(x,p)-\varepsilon \div(p)$ is 
\begin{equation*} 
N(x,p)=g(x)+G(p_1^-,p_2^+,p_3^-,p_4^+)\\-\varepsilon \left(\frac{p_1-p_2}{\Delta x}+\frac{p_3-p_4}{\Delta x}\right),
\end{equation*}
where for a real number $r$, we define the operators
\begin{equation}\label{eq:monotonicity operators}
r^+ :=\max(0,r), \ \ \ r^-:=\max(0,-r),
\end{equation}
and 
$$G(p)=G(p_1,p_2,p_3,p_4):=(p^2_1+p_2^2+p^2_3+p_4^2)^\frac{\alpha}{2}.$$
The operators $r^+$ and $r^-$ are chosen to preserve the monotonicity of the scheme for the HJ operator, which is well defined backward in time.

Now, we compute the operator $K(x,\mathcal D U,M)$, 
and we obtain
\begin{multline*}
K(x_{i,j},[\mathcal D U]_{i,j},M_{i,j})= \\[6pt]
\frac{1}{\Delta x} \left[M_{i,j}\frac{\partial N}{\partial p_1}(x_{i,j},[\mathcal D U]_{i,j})-M_{i-1,j}\frac{\partial N}{\partial p_1}(x_{i-1,j},[\mathcal D U]_{i-1,j})\right.\\[6pt]
+M_{i+1,j}\frac{\partial N}{\partial p_2}(x_{i+1,j},[\mathcal D U]_{i+1,j})-M_{i,j}\frac{\partial N}{\partial p_2}(x_{i,j},[\mathcal D U]_{i,j})\\[6pt]
+M_{i,j}\frac{\partial N}{\partial p_3}(x_{i,j},[\mathcal D U]_{i,j})-M_{i,j-1}\frac{\partial N}{\partial p_3}(x_{i,j-1},[\mathcal D U]_{i,j-1})\\[6pt]
\left.+M_{i,j+1}\frac{\partial N}{\partial p_4}(x_{i,j+1},[\mathcal D U]_{i,j+1})-M_{i,j}\frac{\partial N}{\partial p_4}(x_{i,j},[\mathcal D U]_{i,j})  \right] \\[6pt]
-\varepsilon \ \Delta_d M_{i,j}.
\end{multline*}
We use this operator in \eqref{fp3}. This scheme is similar to the one in \cite{AchdouCapuzzo10}.


%

\subsection{Semi-Lagrangian scheme}
To describe a semi-Lagrangian scheme appropriate to approximate \eqref{ham}, we introduce the operator 
\begin{equation}
\mathcal D^{\gamma}u:=\max_{\gamma\in B(0,1)}\frac{\mathcal I[u](x,\gamma)-u(x)}{\Delta x},
\end{equation}
where $B(0,1)$ is the unitary ball in $\RR^2$, and 
\begin{equation*} 
\mathcal I [u](x,\gamma)=\frac{1}{2}\sum_{i=1}^2 \left( \mathbb I[u](x+\gamma\Delta x +e_i\sqrt{2\varepsilon\Delta x})\right.\\
\left.+\mathbb I[u](x+\gamma\Delta x -e_i\sqrt{2\varepsilon\Delta x})\right).
\end{equation*}
Here, $\mathbb I[u](x)$ is an interpolation operator on the matrix $u$, and $e_i$ is the $i$ unitary vector of an orthonormal basis of the space. The approximation of  $H(x,p)-\varepsilon \div(p)$ is then simply
\begin{equation*} 
N(x,p)=g(x)+p^\alpha.
\end{equation*}
We take the adjoint of the linearized of $N$, 
and we use it into \eqref{fp3}, analogously as performed for the finite-difference scheme.
This scheme differs from the one proposed in \cite{carlini2016DGA}, where an estimation on the volumes of the density distribution $M$ was necessary. We note that the operator $N(x,p)$ is monotone, see \cite{falcone2013semi}.
%
%
\section{Applications to Systems of PDEs} \label{mf}

One immediate application of our numerical scheme is to solve "measure-potential" systems of PDEs. These systems comprise an equation for the evolution of a measure coupled with a second equation for a potential or value function. Typically, this potential determines the drift for the convection in the first equation. Many problems have this structure: mean-field games, traffic-flow models, crowd motion, and chemotaxis. Here, we describe how to use our method in the following examples: two 1-D forward-forward mean-field game (FFMFG) problems and a 2-D crowd motion model.

\subsection{Example: 1-D forward-forward mean-field games}\label{sct:1dmfg}
Here, we consider two one-dimensional forward-forward mean-field game problems, see \cite{AchdouCapuzzo10,MR3575617,GomesSedjroFFMFGCongestion}. The general form of such systems is
\begin{equation}\label{sys:FF MFG}
\begin{cases}
u_t + H(u_x)= \varepsilon u_{xx} + g(\rho), \\[6pt]
\rho_t-(H'(u_x) \rho)_x=\varepsilon \rho_{xx},   \end{cases}
\end{equation}
together with the {\em initial-initial conditions}:  
\begin{equation*}\label{ini-ini}
\begin{cases}
u(x,0) = u_0(x), \\[5pt]
\rho(x,0) = \rho_0 (x).
\end{cases}
\end{equation*}

\begin{figure}[htb]
	\centering
	\begin{subfigure}[b]{\sizefigure\textwidth}
		\includegraphics[width=\textwidth]{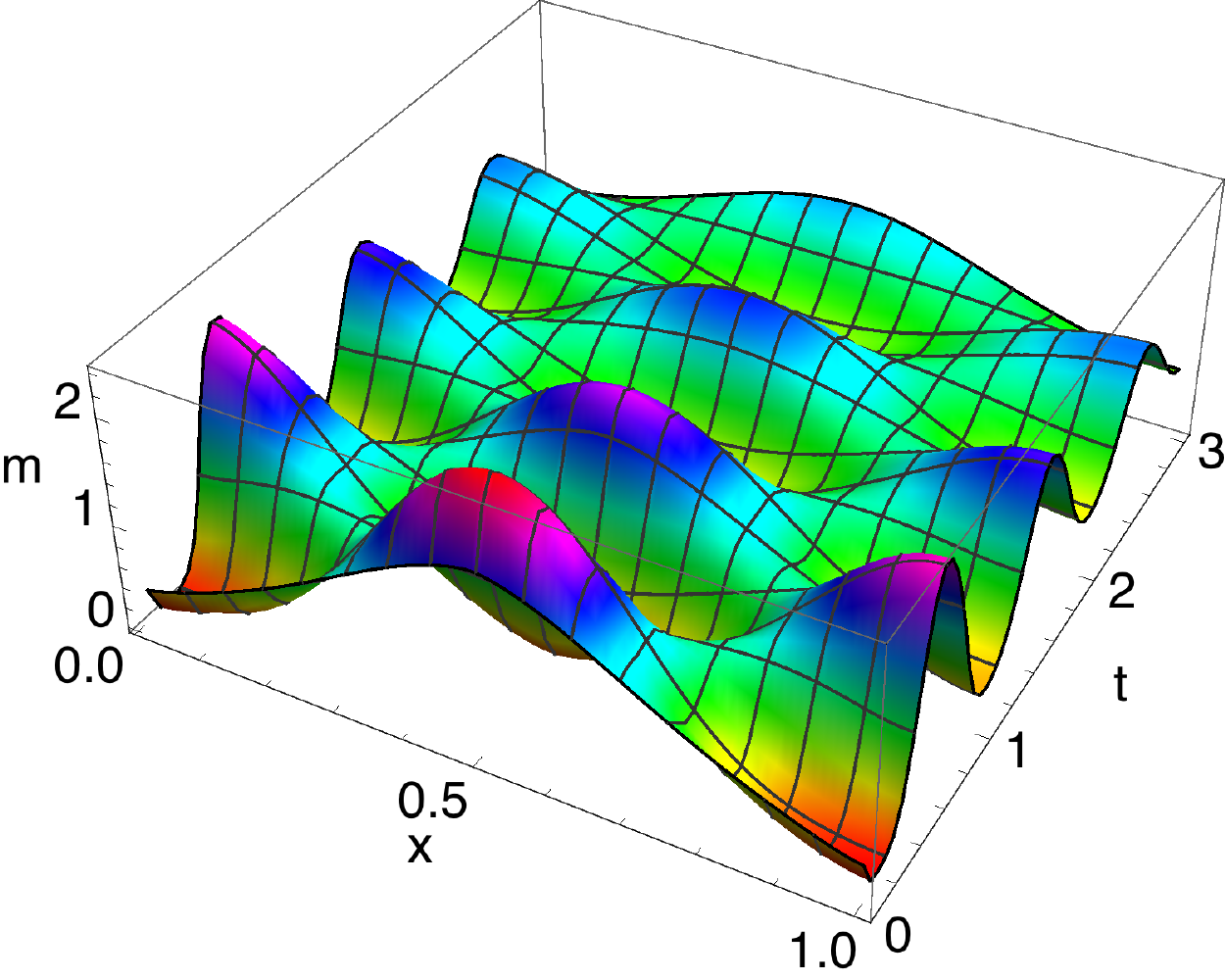}
		\caption{Density}
		\label{fig:solDensity1Drho}
	\end{subfigure}
	~ 
	\begin{subfigure}[b]{\sizefigure\textwidth}
		\includegraphics[width=\textwidth]{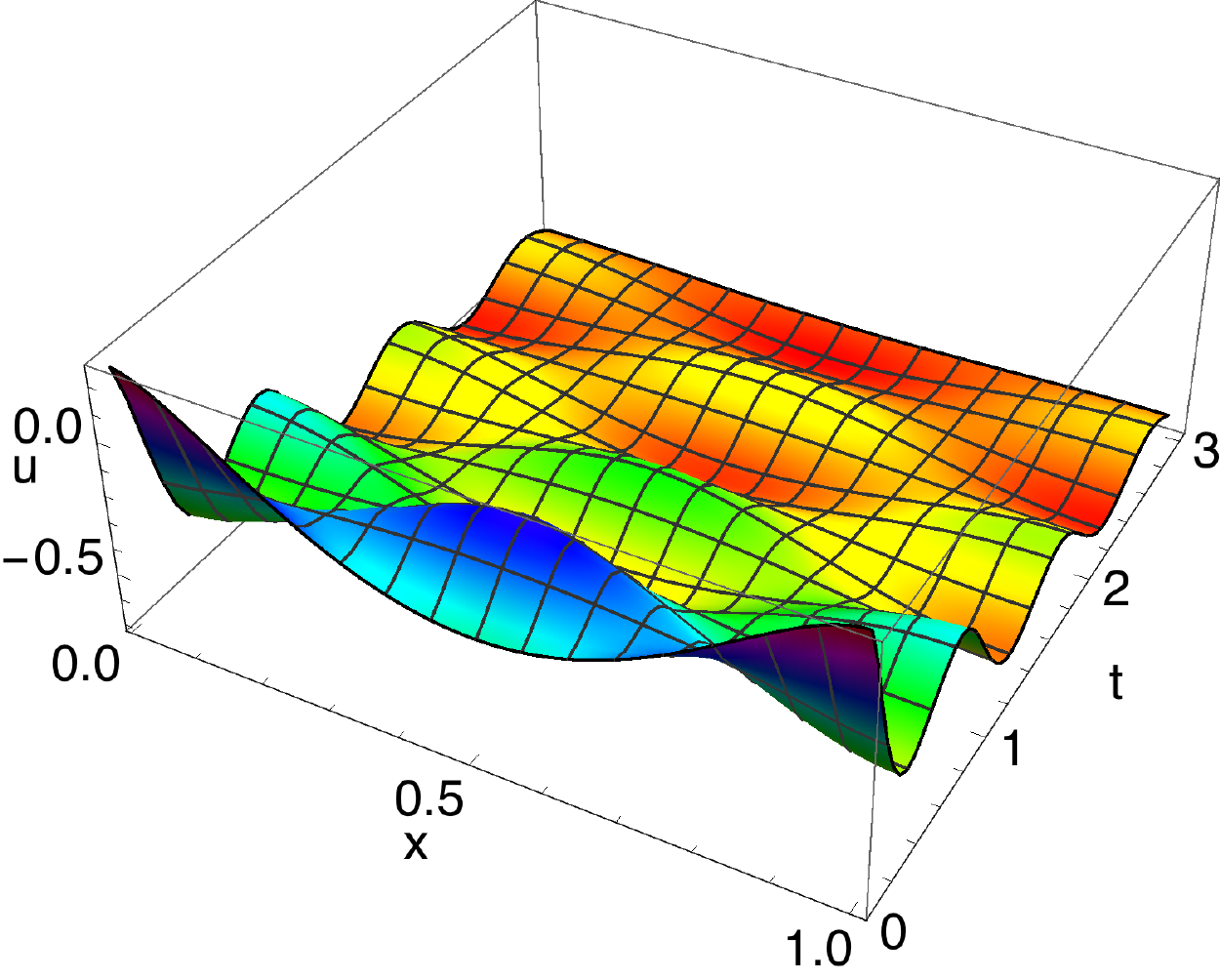}
		\caption{Value function}
		\label{fig:solEikonal1Drho}
	\end{subfigure}
	\caption{Solutions for $g(\rho) = \ln \rho$.}\label{fig:sol FF MFG log rho}
\end{figure}

In this example, we use periodic boundary conditions. For the first problem, we set $H(u_x) =  \frac{u_x^2}{2}$, $g(\rho) = \ln \rho$, and $\varepsilon=0.01$. We then solve:
\begin{equation}\label{sys:FF MFG Solved1}
\begin{cases}
u_t + \frac{u_x^2}{2} = 0.01 \ u_{xx} + \ln \rho,\\[6pt]
\rho_t-(u_x \rho)_x= 0.01 \ \rho_{xx}.   \end{cases}
\end{equation}
We choose the {\em initial-initial conditions}:
\begin{equation*}
\begin{cases}
u_0(x)=0.3 \cos(2\pi x), \\[5pt]
\rho_0(x)=1.
\end{cases}
\end{equation*}
We depict the solution of this problem in Figure~\ref{fig:sol FF MFG log rho}. 

%
%
%


\begin{figure}[htb]
	\centering
	\begin{subfigure}[b]{\sizefigure\textwidth}
		\includegraphics[width=\textwidth]{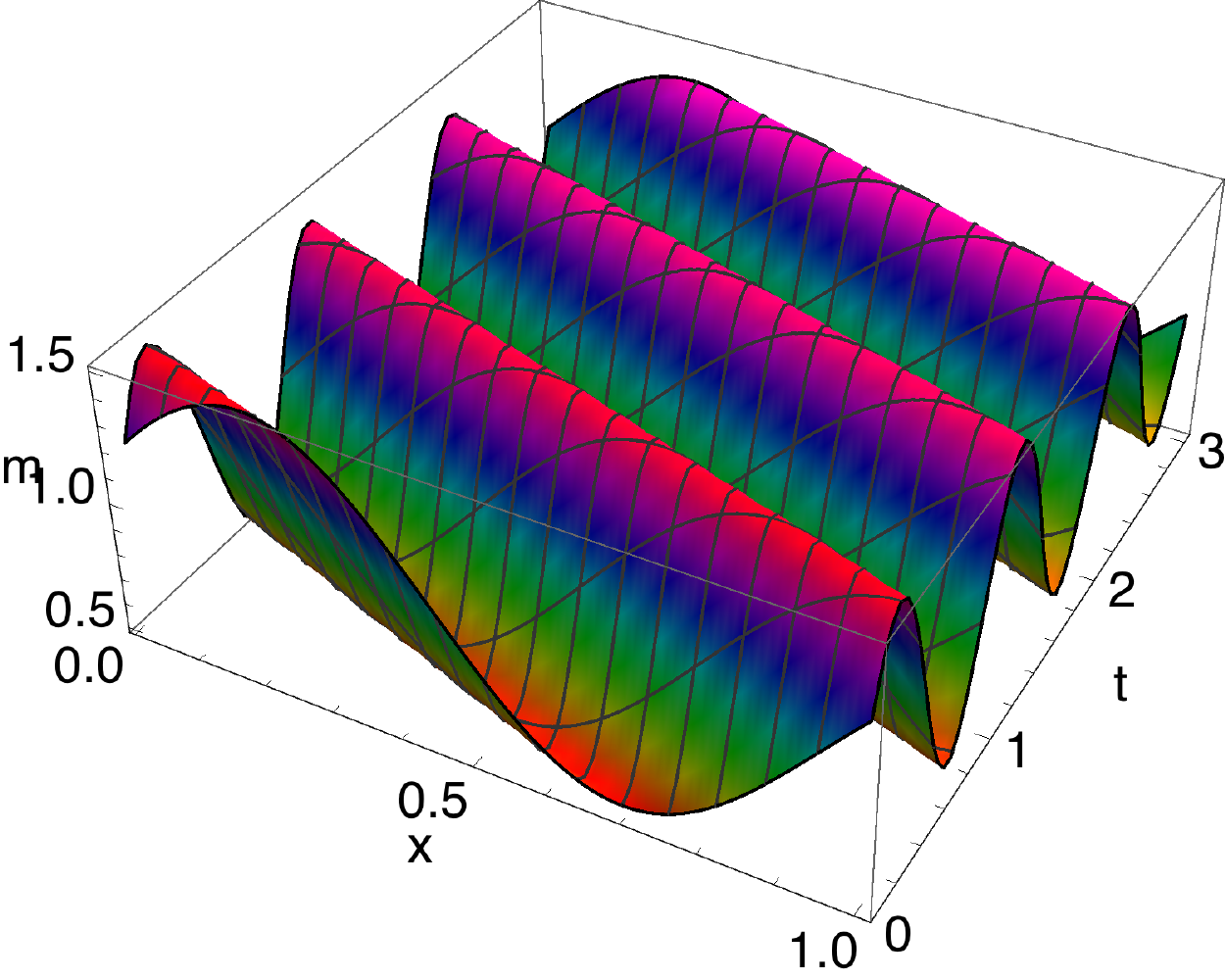}
		\caption{Density}
	\end{subfigure}
	~ 
	\begin{subfigure}[b]{\sizefigure\textwidth}
		\includegraphics[width=\textwidth]{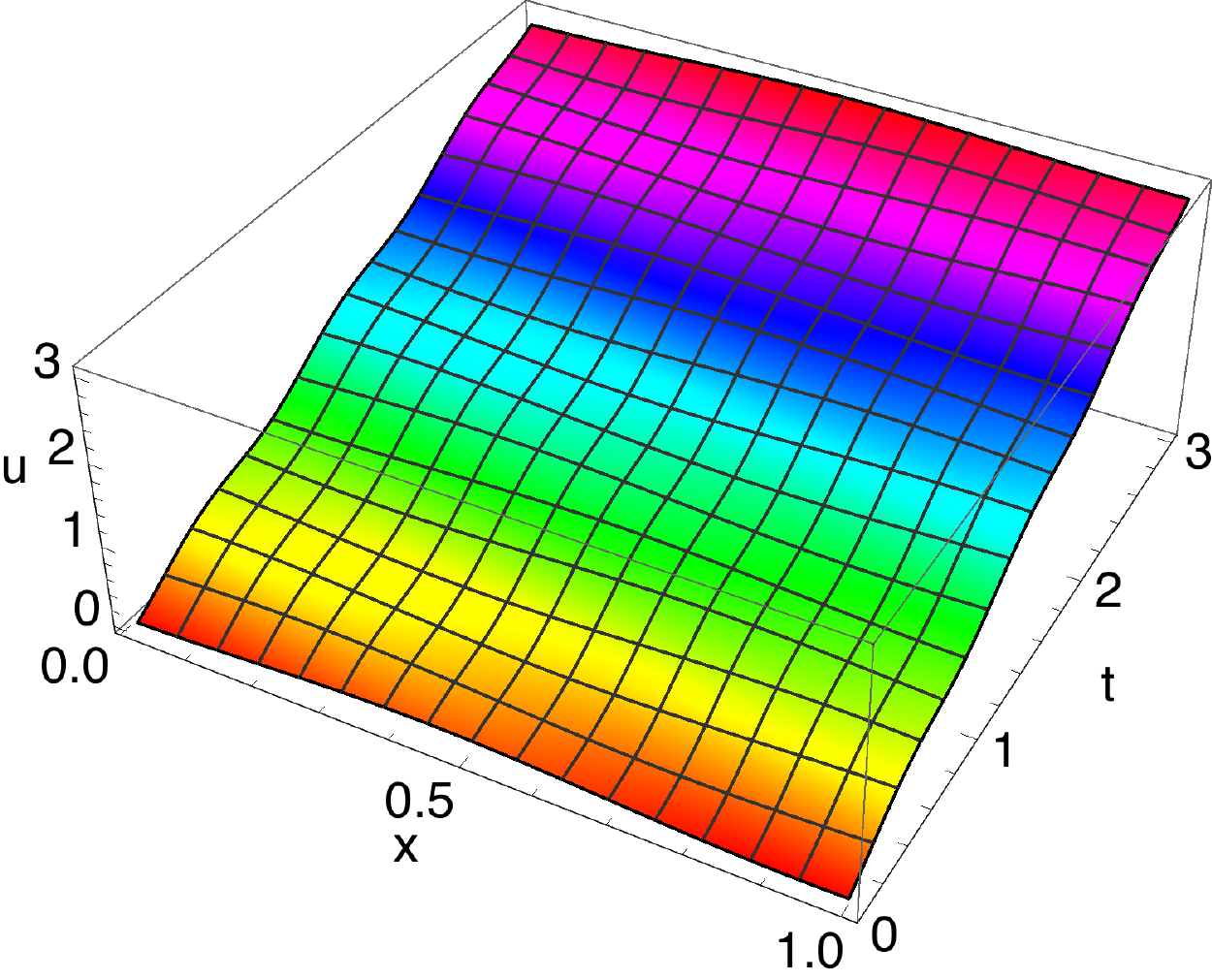}
		\caption{Value function}
	\end{subfigure}
	\caption{Solutions for the FFMFG with congestion.}\label{fig:sol FF MFG Travelling Wave}
\end{figure}

Now, for the second case, we choose $H(u_x,\rho) = \frac{(p+u_x)^2}{2\rho^\alpha}$, $g(\rho)=\frac 3 2 \rho^\alpha$, and $\varepsilon=0$. This is a first-order FFMFG with congestion, which is equivalent to a system of conservation laws. Setting $v = p + u_x$, the equivalent system is
\begin{equation}\label{sys:FF MFG Solved2V}
\begin{cases}
v_t + \left(\frac{v^2}{2 \rho^\alpha} - \frac{3}{2} \rho^\alpha\right)_x = 0, \\[6pt]
\rho_t- \left( \rho^{1-\alpha}v \right)_x=0.   \end{cases}
\end{equation}
For $\alpha=1$, and for the {\em initial-initial conditions}
\begin{equation*}
\left\{\begin{array}{ll}
u_0= -0.5 \ \frac{\cos(2 \pi x)}{2 \pi},& \\[5pt]
\rho_0=1+0.5 \sin(2 \pi x),&
\end{array}\right.
\end{equation*}
the solution for the density in \eqref{sys:FF MFG Solved2V} is a traveling wave; as shown in \cite{GomesSedjroFFMFGCongestion}, and depicted in Figure~\ref{fig:sol FF MFG Travelling Wave}.

%
%
%
%

Now, we explain how we treated such systems numerically. MFGs have built-in the adjoint structure we consider here. Hence,  we can use the same spatial discretization for both the FP and HJ equations. Each of the discretizations requires solving an ODE in time. Since we must solve the system of FP coupled to a HJ equation, we treat these ODEs as a system, and we can apply a suitable solver for the time discretization. In our examples, we use finite differences for the spatial discretization,  as in section~\ref{ssct:FinDif}. The simulations corresponding to 
Figure \ref{fig:sol FF MFG log rho} and Figure~\ref{fig:sol FF MFG Travelling Wave} were produced with a spatial grid with 80 points, final time $T~=~3$, and 50 points for the sample on time.


\subsection{Example: Hughes Model in 2-D}\label{Example_Hughes}

In this example, we present a model for crowd motion model due to Hughes \cite{hughes2002continuum,hughes2000flow}. The model comprises a FP equation, describing the evolution of the density of pedestrians/agents, coupled to an Eikonal (EK) equation that gives the optimal movement direction. This two-dimensional system is
\begin{equation}\label{Hughes_System}
\begin{cases}
\displaystyle \rho_t(x,t) - \div(\rho(1-\rho)^2 Du)= 0,\\[6pt] 
\displaystyle |Du(x)|^2=\frac{1}{(1-\rho)^2},  
\end{cases}
\end{equation} 
together with an initial condition for the density. The goal is to exit a domain $\Omega$ in minimal time taking into account congestion effects. Due to the stationary character of the EK equation, this system is not of mean-field game type. The density, $\rho$, evolves as if at each instant of time the EK equation sees a frozen density. Then the agents choose the direction that leads to the shortest-time to evacuation and this process determines the evolution of $\rho$. 

Now, we describe how the Hughes system fits our framework. Performing the same steps as in section \ref{eaa}, with the HJ operator
\begin{equation}\label{HJ_with_f}
- u_t + f(\rho) H(x,Du) - \varepsilon \Delta u,
\end{equation}
where $f(\rho)$ is a regular function of the density, we obtain the associated FP equation
\begin{equation}\label{FP_with_f}
 \rho_t - \div \left( f(\rho) D_pH(x,Du) \rho \right) = \varepsilon \Delta u.
\end{equation}
By setting $f(\rho) = (1-\rho)^2$ and $H(x,p) = \displaystyle \frac{|p|^2}{2}$, \eqref{FP_with_f} becomes the first equation of \eqref{Hughes_System}; and \eqref{HJ_with_f} is the adjoint operator we must study. Since the EK equation is a particular case of a HJ equation, we discretize it in space as with the HJ operator associated to the FP equation. In the following example, we use finite differences to discretize the generator of the HJ operator. For the time discretization, we use an explicit Euler method. 

The domain is a rectangle $[0,3] \times [0,1]$, with an exit on $[2.25,3] \times \{1\}$, corresponding to a typical proportional size of a door in a room. We set the value of $u$ to $+\infty$ on all the boundary but on its exit, where we fix it equal to zero. The density is set equal zero on the boundary.

\begin{figure}[htb]
	\centering
	\begin{subfigure}[b]{\sizefigure\textwidth}
		\includegraphics[width=\textwidth]{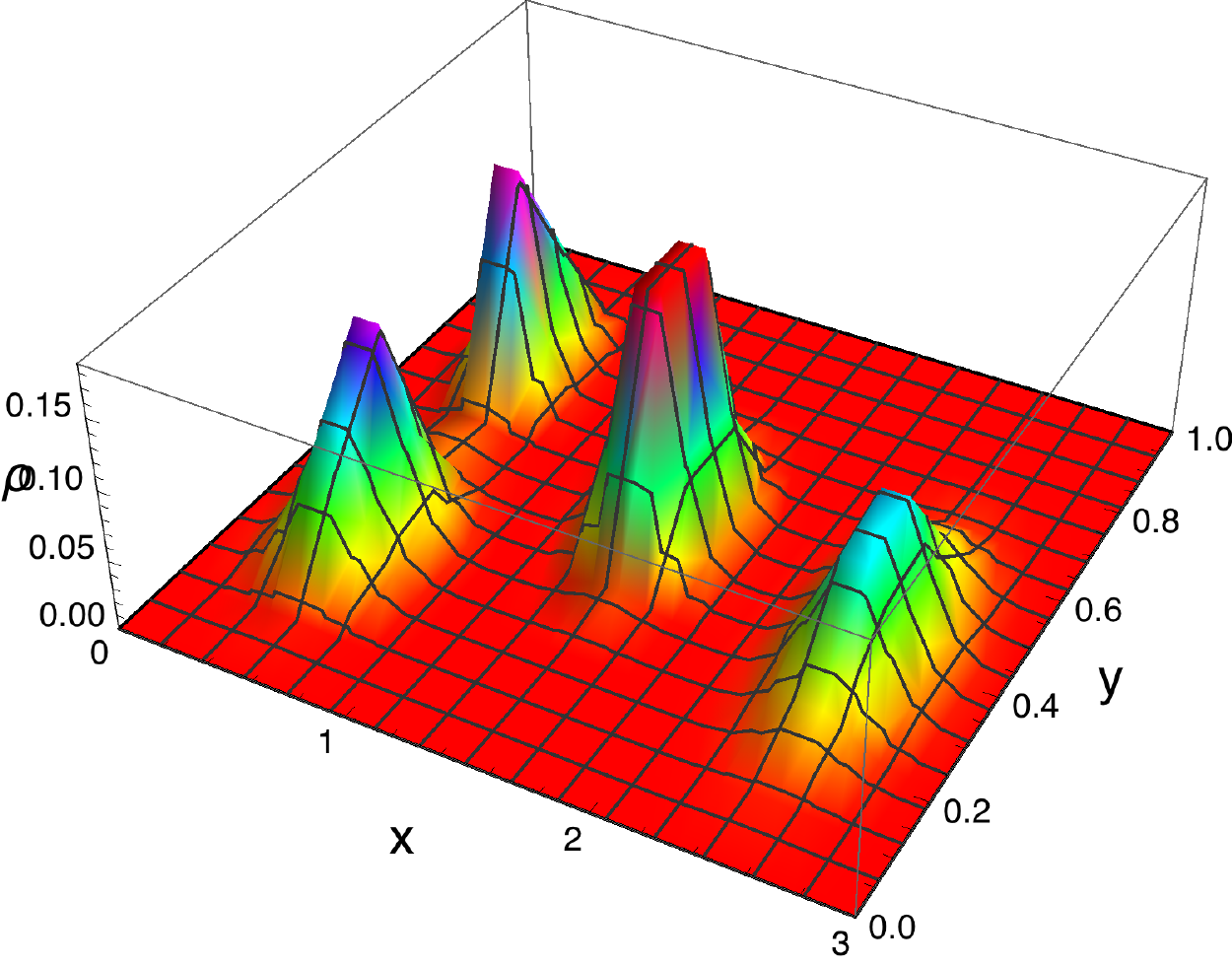}
		\caption{Initial Density.}
	\end{subfigure}
	~ 
	\begin{subfigure}[b]{\sizefigure\textwidth}
		\includegraphics[width=\textwidth]{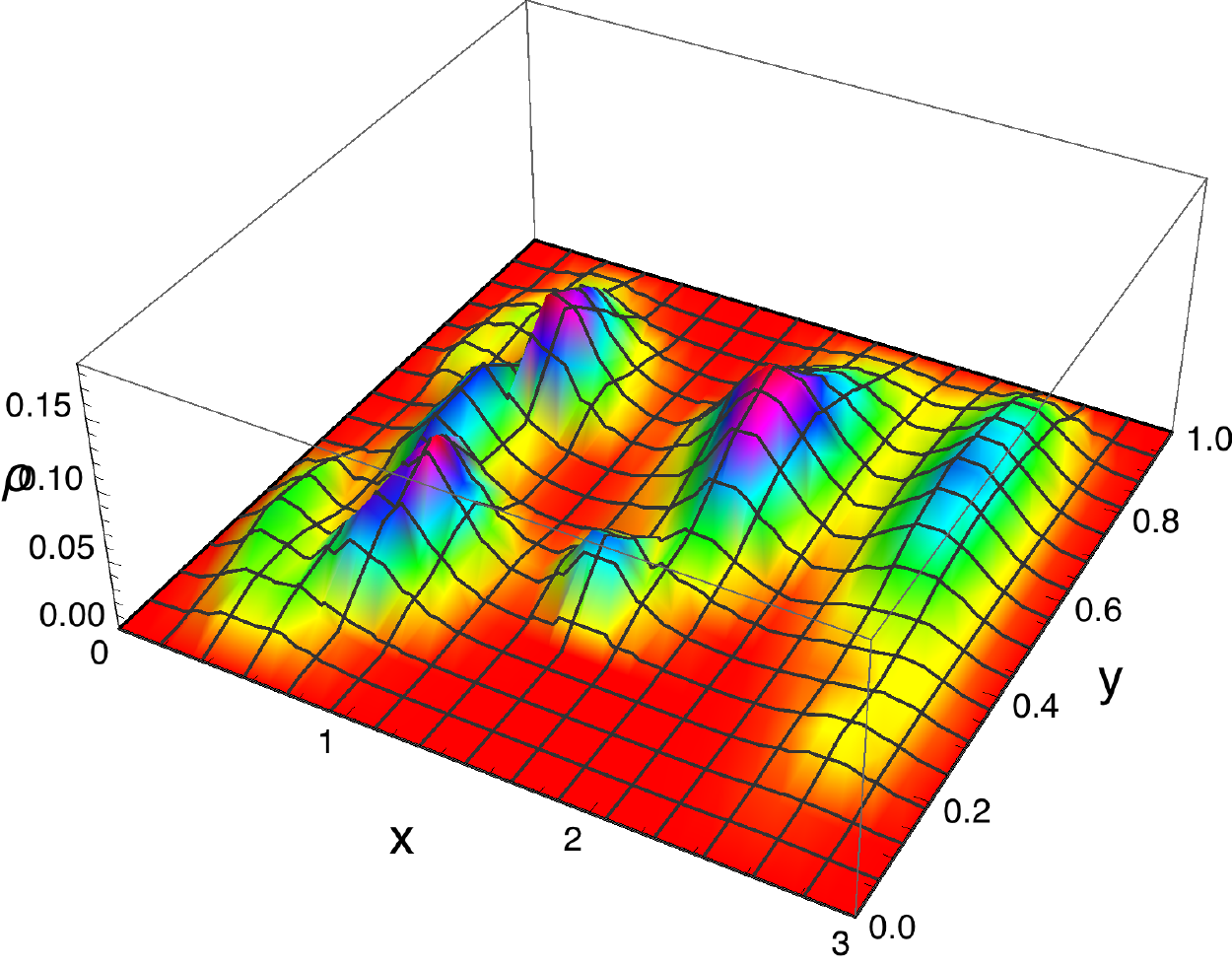}
		\caption{Density at time $0.33$.}
	\end{subfigure} 
\qquad
	\begin{subfigure}[b]{\sizefigure\textwidth}
		\includegraphics[width=\textwidth]{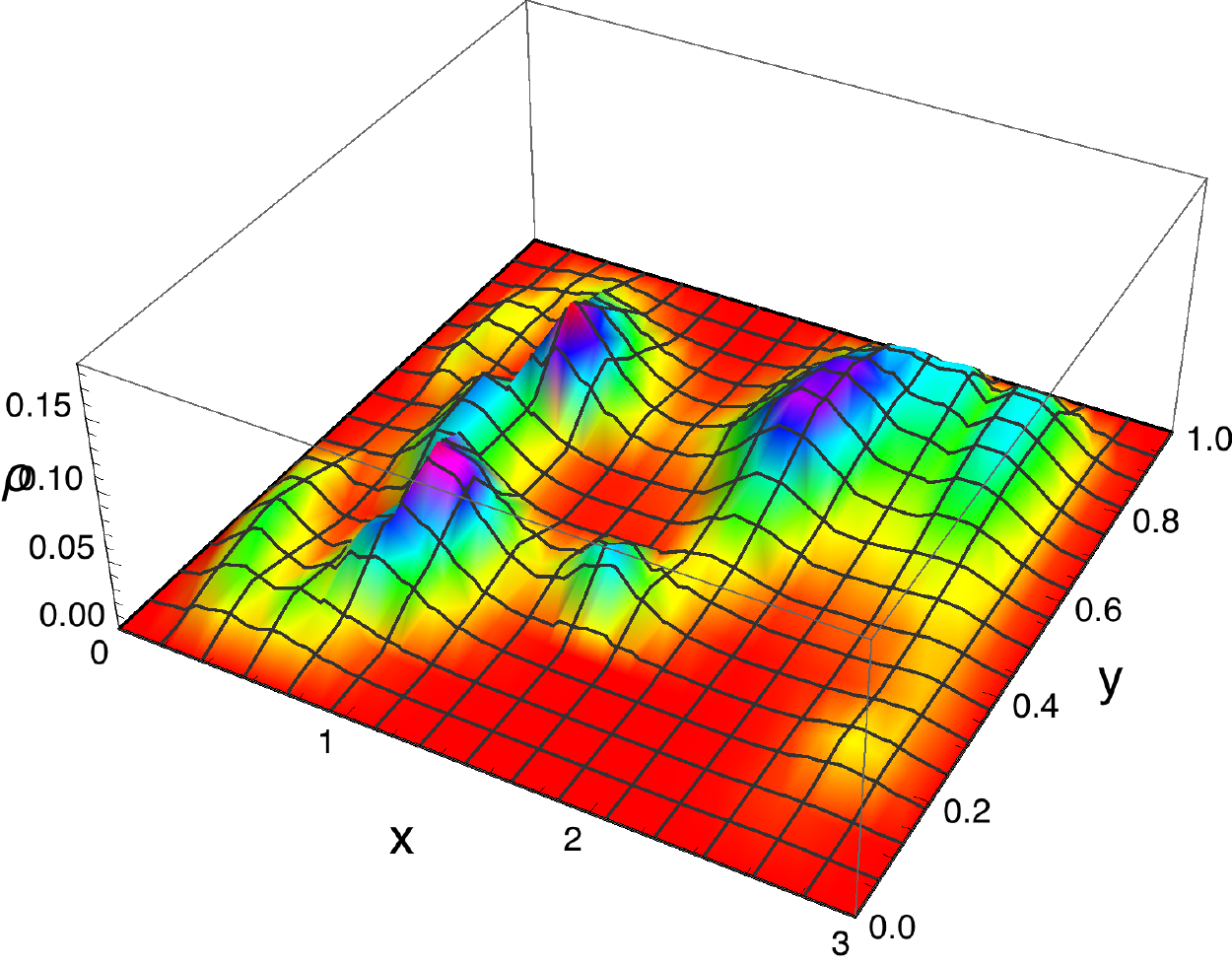}
		\caption{Density at time $0.5$.}
	\end{subfigure}
	~
	\begin{subfigure}[b]{\sizefigure\textwidth}
		\includegraphics[width=\textwidth,scale=0.5]{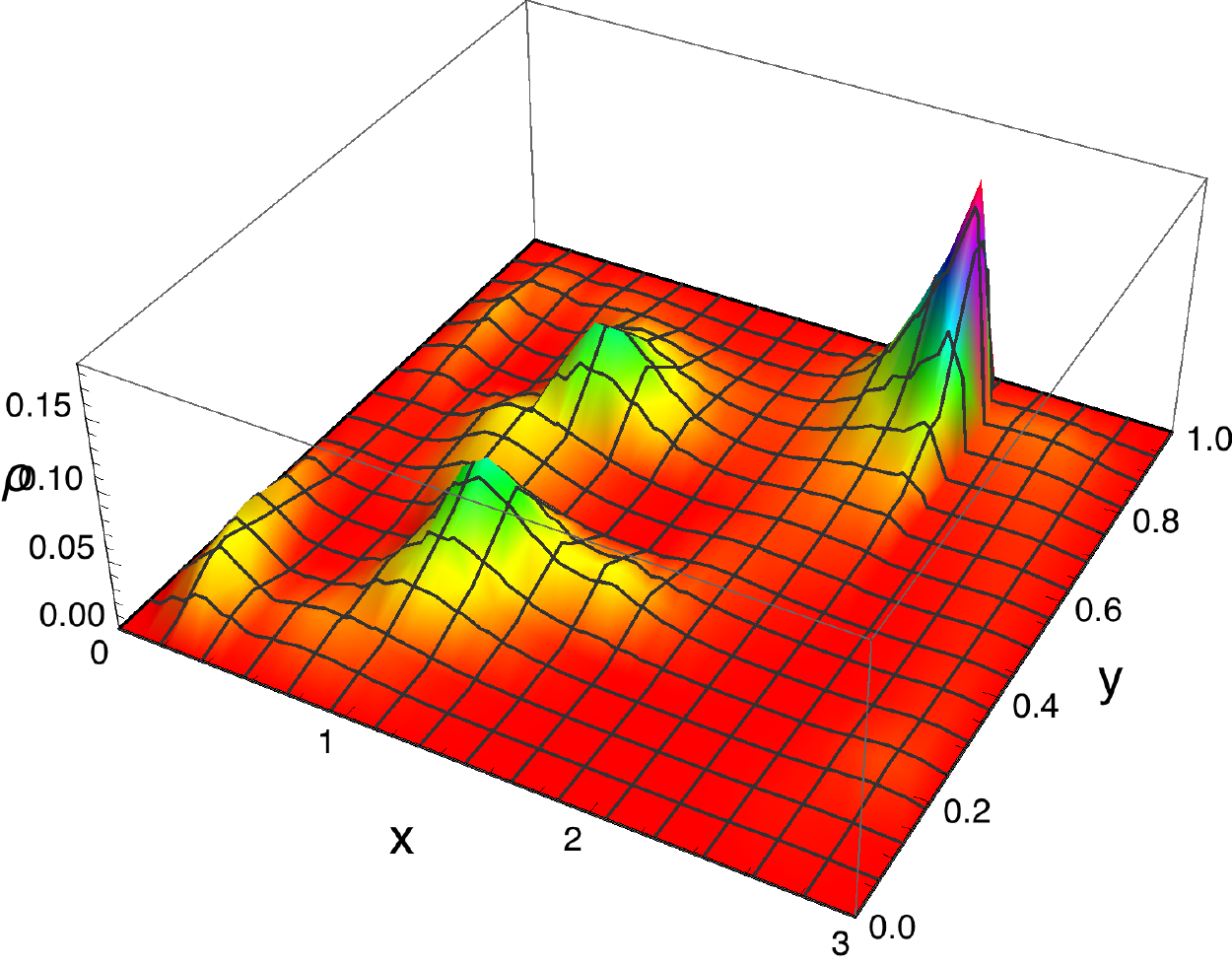}
		\caption{Final density at time $1.0$.}
	\end{subfigure}
	\caption{Evolution of the density for the Hughes model.} \label{fig:Hughes_2D}
\end{figure}

In contrast with MFG problems, the Hughes model does not have the adjoint structure built-in. Again, the numerical solution of the FP equation requires solving an ODE in time. However, the EK equation must be treated in another way; at each iteration of the solver for the FP equation, we solve the EK equation. We use a fixed-point approach, as described in \cite{MR2218974}. Alternatively, fast marching or policy iteration methods could also be applied. We depict the initial condition and its evolution in Figure~\ref{fig:Hughes_2D}. The spatial grid contains $100$ points, and we choose the final time $T=1.0$. 

We end this section by remarking that in the last three problems our simulations preserve mass and positivity, as expected.

%
%
%
%
%
%

\section{Conclusions}
Here, we develop numerical methods to solve nonlinear Fokker-Planck equations via its adjoint Hamilton-Jacobi operator. Our method preserves mass and positivity, and we use it to solve systems of PDEs with a Fokker-Planck equation coupled to a Hamilton-Jacobi equation. Our methods apply to a broad range of problems with a measure-potential structure that include mean-field games, crowd and traffic models, and chemotaxis. 
  
In future work, we plan to 
address different schemes developed for HJ equations to study FP equations. Thus, reversing the process that gave rise to effective numerical schemes for HJ equations, as Discontinuous Galerkin or ENO schemes, originally developed for conservation laws. Nevertheless, it is clear that, without monotonicity and stability properties, results for the convergence of such schemes are difficult to achieve.

\bibliography{mfg}
\end{document}